\documentclass[oneside,english]{article} 
\usepackage[utf8]{inputenc} 
\usepackage{amssymb,amsthm,amsmath} 
\usepackage{babel}
\usepackage{natbib}

\theoremstyle{plain} 
\newtheorem{thm}{Theorem}[section] \numberwithin{equation}{section} \numberwithin{figure}{section} \theoremstyle{plain}

\newtheorem{fact}[thm]{Fact} \theoremstyle{plain} 
\newtheorem{prep}[thm]{Proposition} \theoremstyle{plain} 
\newtheorem{definition}[thm]{Definition} \theoremstyle{remark}
\newtheorem{rem}[thm]{Remark} \theoremstyle{pain}  
 \theoremstyle{remark}

\newtheorem*{acknowledgement*}{Acknowledgement}

\newcommand{\intr}{\int_{\mathbb{R}^d}} 
\newcommand{\inti}{\int_{0}^{+\infty}} 
\newcommand{\intc}[1]{\int_{0}^{#1}} 
\newcommand{\T}[1]{\mathcal{T}_{#1}} 
\newcommand{\pb}{{1+\beta}} 
\newcommand{\norm}[3]{\Vert #1 \Vert_{ #2 } ^{ #3 }} 
 
\newcommand{\rbr}[1]{\left( #1 \right)}

\newcommand{\sgn}{\text{sgn}}

\newcommand{\ddp}[2]{\left\langle #1, #2 \right\rangle}

\newcommand{\Rd}{\mathbb{R}^d}

\author{Piotr Mi\l{}o\'s}
\newcommand{\ev}[1]{\mathbb{E}{#1}} 
\title{Occupation time fluctuation limits of infinite variance equilibrium branching systems}

\begin{document}
\maketitle
	\begin{abstract}
		We establish limit theorems for the fluctuations of the rescaled occupation time of a $(d,\alpha,\beta)$-branching particle system. It consists of particles moving according to a symmetric $\alpha$-stable motion in $\mathbb{R}^d$. The branching law is in the domain of attraction of a (1+$\beta$)-stable law and the initial condition is an equilibrium random measure for the system (defined below). In the paper we treat separately the cases of intermediate $\alpha/\beta<d<(1+\beta)\alpha/\beta$, critical $d=(1+\beta)\alpha/\beta$ and large $d>(1+\beta)\alpha/\beta $ dimensions. In the most interesting case of intermediate dimensions we obtain a version of a fractional stable motion. The long-range dependence structure of this process is also studied. Contrary to this case, limit processes in critical and large dimensions have independent increments.
	\end{abstract}
	
	
	
	\section{Introduction} 
	\subsection{Branching system and occupation time fluctuations} 
	The aim of this paper is to present some (functional) limit theorems for the occupation time fluctuation process of a branching particle system. We call a $(d, \alpha, \beta)$-branching particle system (denoted in the sequel by $N$) a set of particles moving independently according to the spherically symmetric $\alpha$-stable L\'evy motion ($0<\alpha\leq 2$) in $\mathbb{R}^d$ and splitting after exponential time (with intensity $V$) with branching law
	\[ p_k=\left\{ 
	\begin{array}{cc}
		0 & k=1\\
		\frac{1}{1+\beta} \binom{\pb}{k}(-1)^k & k = 0,2,3,\ldots 
	\end{array}
	\right. \]
	$( 0<\beta < 1 )$. This is an example of a law in the domain of attraction of $(1+\beta)$-stable variable. It has infinite variance and is critical. For $\beta=1$ it reduces to binary critical branching which was treated in a series of papers mentioned below. The generating function of this law is
	\begin{equation}
	F(s) = s + \frac{1}{1+\beta} (1-s)^{\pb}, s\in(0,1). \label{def:generating_beta}
	\end{equation}
The particle system will be represented by an empirical measure process $(N_t)_{t \geq 0}$, i.e. for a Borel set $A$, $N_t(A)$ is a (random) number of particles in $A$ at time $t$.
	The initial particle distribution is yet to be introduced. The most natural choice is a Poisson random field with homogeneous intensity, i.e. (Lebesgue measure) $\lambda$. This case, which was studied in \cite{Bojdecki:2007aa} and \cite{bojdecki-2005}, is a starting point and reference for our investigation. It is known \cite{Gorostiza:1991aa} that for $\alpha/\beta<d$ such system (denoted by $N^{Poiss}$) converges to an equilibrium distribution
	\begin{equation}
	 N^{Poiss} \Rightarrow Eq \label{def:Eq}
	\end{equation}
	where $\Rightarrow$ denotes weak convergence in the space of point measures. The Laplace functional of the equilibrium distribution is given by
\begin{equation}
	\mathbb{E}\exp\left\{ -\left\langle Eq,\varphi \right\rangle \right\} =\exp\left\{ \left\langle \lambda,e^{-\varphi}-1\right\rangle +V\int_{0}^{\infty}\left\langle \lambda,H\left(j\left(\cdot,s\right)\right)\right\rangle ds\right\} ,\label{eq: laplace_equilibrum}\end{equation}
where 
\begin{equation}
	j\left(x,l\right):=\mathbb{E}\exp\left(-\left\langle N_{l}^{x},\varphi \right\rangle \right)
	\label{def: h}
\end{equation}
$N^x$ is the empirical process of the system starting from $x\in \Rd$, 
 $H(s)=F\left(s\right)-s$, $\varphi:\mathbb{R}^d\rightarrow\mathbb{R}_{+}$, $\varphi \in \mathcal{L}^1(\mathbb{R}^d)\cap C(\mathbb{R}^d)$ and $j$ satisfies the integral equation\[
j\left(x,l\right)=\mathcal{T}_{l}e^{-\varphi}\left(x\right)+V\int_{0}^{l}\mathcal{T}_{l-s}H\left(j\left(\cdot,s\right)\right)\left(x\right)ds.\]
 This equations can be obtained in the same way as \cite[(2.4)]{Gorostiza:1991aa}. 
In this paper we consider a system $N$ starting off from $Eq$ and compare the obtained result to the ones in \cite{Bojdecki:2007aa} and \cite{bojdecki-2005}. 
	 For the process $(N_t)_{t\geq 0}$ we define the rescaled occupation time fluctuations process by
	\begin{equation}
	 X_T(t) = \frac{1}{F_T} \int_0^{Tt} (N_s - \lambda) ds, \label{def:occupation_fluct}
	\end{equation}
	where $F_T$ is a proper normalization and $T$ is a scaling parameter which accelerates the time. The object of our investigation is the limit of  $X_T$ as $T$ tends to $+\infty$
\begin{equation}
	X_T \Rightarrow X. \label{lim:main}
\end{equation}
	For the time being we are not very rigorous and do not specify the type of convergence.
	\subsection{Results and proof techniques}
	In the proofs we will rely on methods presented in \cite{bojdecki-2005}, \cite{Bojdecki:2007aa} and \cite{Milos:2007aa}.\\
	 Although the process $X_T$ is signed-measure-valued it is convenient to regard it as a process with values in the space $\mathcal{S}'(\Rd)$ of tempered distributions, which is dual to the space of smooth and rapidly decreasing functions $\mathcal{S}(\Rd)$. We denote duality in this space by $\ddp{\cdot}{\cdot}$. In this space one may employ space-time method introduced by \cite{Bojdecki:1986aa} which together with Mitoma's theorem constitute a powerful technique in proving weak, functional convergence.\\
	Three kinds of convergence are used. The convergence of finite-dimensional distributions is denoted by $\Rightarrow_f$. For a continuous, $S'(\Rd)$-valued process $X = (X_t)_{t\geq 0}$ and any $\tau>0$ one can define an $\mathcal{S}'(\mathbb{R}^{d+1})$-valued random variable 
	\begin{equation}
		\ddp{\tilde{X}}{\Phi} = \intc{\tau} \ddp{X_s}{\Phi(\cdot, s)} ds, \: \Phi \in \mathcal{S}(\mathbb{R}^{d+1}). \label{def: space-time} 
	\end{equation}
	If for any $\tau>0$ $\tilde{X}_n\rightarrow \tilde{X}$ in distribution, we say that the convergence in the space-time sense holds and denote this fact by $\Rightarrow_i$. Finally, we consider the functional weak convergence denoted by $X_n\Rightarrow_c X$. It holds if for any $\tau>0$ processes $X_n = (X_n(t))_{t\in [0, \tau]}$ converge to $X = (X(t))_{t\in [0, \tau]}$ weakly in $C([0, \tau], \mathcal{S}'(\Rd))$ (in the sequel without loss of generality we assume $\tau = 1$). It is known that $\Rightarrow_i$ and $\Rightarrow_f$ do not imply each other, but either of them together with tightness implies $\Rightarrow_c$ \cite{Bojdecki:1986aa}. Conversely, $\Rightarrow_c$ implies both $\Rightarrow_i$, $\Rightarrow_f$. \\
	The presentation of the results naturally splits into parts, corresponding to intermediate dimensions
	\[ \frac{\alpha}{\beta} < d < \frac{\alpha (\pb)}{\beta}, \]
 critical
	\[ d = \frac{\alpha (\pb)}{\beta}, \]
	and large dimensions
	\[ d > \frac{\alpha (\pb)}{\beta}, \]
	respectively.\\
	In the first case of intermediate dimensions we obtain a weak functional convergence to process $X$ of the form $X=K\eta\lambda$ where $K$ is a constant, $\eta$-a stable process being in a sense a stable (non-Gaussian) analogue of a fractional Brownian motion. 
So we see that in this case the limit has a very simple spatial structure whereas its temporal structure is complicated.
It is also worthwhile to point out that $\eta$ has a stationary increments (unlike the corresponding process in \cite{Bojdecki:2007aa}) and is a heavy-tailed process with long-range dependence. This dependence is described in Section \ref{sec:results} in terms of dependence exponent  and roughly speaking it means that "dependence" decays polynomially.
	The cases of critical and large dimensions differ substantially, one can prove only finite-dimensional distributions and space-time convergences. In both cases we obtain processes with independent increments and the limit for large dimensions is truly $\mathcal{S}'(\Rd)$-valued. These processes are not continuous hence it is not possible to obtain functional convergence.
	\subsection{A survey of results}

	This paper is a part of a larger programme carried out by Bojdecki et al. and recently by Milos. It seems useful to present it here in a compact way. This small survey is not meant to be exhaustive nor very strict it aims only to give a reader a glimpse of the whole picture. In tables below we gather limits of occupation time fluctuations under time rescaling (as defined by (\ref{lim:main})), type of convergence and the normalizing factor in different settings. The structure of tables reflects dependence on the dimension of the space and the starting distribution.\\
\begin{center}
	\begin{tabular}{|c|c|c|}
      \multicolumn{3}{c}{Table 1. Systems with finite variance branching law}\\
		\hline
		&
		Poisson
		&
		Equilibrium
		\tabularnewline
		\hline
		$\begin{array}{cc} \alpha <d <2\alpha\\ \textit{intermediate}\end{array}$
		&
		$\begin{array}{ccc} K\cdot\textit{sub-frac-BM}\cdot\lambda \\ \text{functional } \: T^{\frac{3-d/\alpha}{2}} 
								\\ \text{\cite{Bojdecki:2006ab} and \cite{Milos:2007aa}}  \end{array}$
		&
		$\begin{array}{ccc} K\cdot\textit{frac-BM}\cdot\lambda \\ \text{functional } \: T^{\frac{3-d/\alpha}{2}} 
								\\ \text{\cite{Milos:2007aa}}  \end{array}$
		\tabularnewline
		\hline
		$\begin{array}{cc} d = 2\alpha\\ \textit{critical}\end{array}$&
		$\begin{array}{ccc} K\cdot\textit{BM}\cdot\lambda \\ \text{functional } \: (T \log T)^{\frac{1}{2}}
								\\ \text{\cite{Bojdecki:2006aa} and \cite{Milos:2007ab}}  \end{array}$
		 &
		$\begin{array}{cc} \textit{the same}
								\\ \text{\cite{Milos:2007ab}}  \end{array}$
		\tabularnewline
		\hline
		$\begin{array}{cc} d > 2\alpha\\ \textit{large}\end{array}$&
		$\begin{array}{ccc} \mathcal{S}'(\Rd)\textit{-BM} \\ \text{functional } \: T^{\frac{1}{2}}
								\\ \text{\cite{Bojdecki:2006aa} and \cite{Milos:2007ab}}  \end{array}$
		 &
		$\begin{array}{cc} \textit{the same*}
								\\ \text{\cite{Milos:2007ab}}  \end{array}$
		\tabularnewline
		\hline		
	\end{tabular}
   
  \end{center}
* - due to technical difficulties functional convergence proved only for finite fourth-moment branching law.\\
$K$ denotes a generic constant \\ 
\textit{BM} - standard real-valued Brownian motion\\
\textit{sub-frac-BM} - sub-fractional Brownian motion (ie. centered Gaussian process with covariance function $s^h + t^h -\frac{1}{2}[(s+t)^h + |s-t|^h]$) \\
\textit{frac-BM} - fractional Brownian motion (ie. centered Gaussian process with covariance function $\frac{1}{2}[s^h + t^h + |s-t|^h]$)\\
$\mathcal{S}'(\Rd)\textit{-BM}$ - centered Gaussian $\mathcal{S}'(\Rd)$-valued process with covariance functional 
\begin{equation}
	Cov\left( \ddp{X_s}{\varphi_1}, \ddp{X_t}{\varphi_2} \right) = 
(s\wedge t) \frac{1}{2\pi} 
\intr \left( \frac{2}{|z|^\alpha} + \frac{Vm}{2|z|^{2\alpha}} \right) \widehat{\varphi_1}(z) \overline{\widehat{\varphi_2}(z)} dz \nonumber,
\end{equation}
where$\ \varphi_1,\varphi_2 \in\mathcal{S}\left(\mathbb{R}^{d}\right)$ and $m$ depends on branching law.
Papers \cite{Bojdecki:2006ab} and \cite{Bojdecki:2006aa} contain also results for systems without branching.

\begin{center}
	\begin{tabular}{|c|c|c|}
\multicolumn{3}{c}{Table 2. Systems with infinite variance branching law - generating function (\ref{def:generating_beta})}\\
		\hline
		&
		Poisson
		&
		Equilibrium
		\tabularnewline
		\hline
		$\begin{array}{cc} \frac{\alpha}{\beta} <d < \frac{\alpha(\pb)}{\beta}\\ \textit{intermediate}\end{array}$
		&
		$\begin{array}{cccc} K\cdot\textit{sub-frac-SM}\cdot\lambda \\ \text{functional } \\ F_{T} = T^{(2+\beta - \frac{d}{\alpha}\beta)/(\pb)}
								\\ \text{\cite{Bojdecki:2007aa}}  \end{array}$
		&
		$\begin{array}{cccc} K\cdot\textit{frac-SM}\cdot\lambda \\ \text{functional } \\ F_{T} = T^{(2+\beta - \frac{d}{\alpha}\beta)/(\pb)}
								\\ \text{this paper}  \end{array}$
		\tabularnewline
		\hline
		$\begin{array}{cc} d = \frac{\alpha(\pb)}{\beta}\\ \textit{critical}\end{array}$&
		$\begin{array}{cccc} K\cdot\textit{SM}\cdot\lambda \\ \text{fin-dims and space-time } \\ (T \log T)^{\frac{1}{\pb}}
								\\ \text{\cite{bojdecki-2005}}  \end{array}$
		 &
		$\begin{array}{cc} \textit{the same}
								\\ \text{this paper}  \end{array}$
		\tabularnewline
		\hline
		$\begin{array}{cc} d > \frac{\alpha(\pb)}{\beta} \\ \textit{large}\end{array}$&
		$\begin{array}{cccc} \mathcal{S}'(\Rd)\textit{-SM} \\ \text{fin-dims and space-time } \\ T^{\frac{1}{\pb}}
								\\ \text{\cite{bojdecki-2005}}  \end{array}$
		 &
		$\begin{array}{cc} \textit{the same}
								\\ \text{this paper}  \end{array}$
		\tabularnewline
		\hline		
	\end{tabular}
\end{center}
Here,\\
 \textit{sub-frac-SM} - "sub-fractional" stable motion, defined by (\ref{def:eta1})\\
\textit{frac-SM} - "fractional" stable motion, defined by (\ref{def:eta})\\
\textit{SM} - stable motion with independent increments, with finite dimensional distributions given by (\ref{def:lim_critical})\\
$\mathcal{S}'(\Rd)$\textit{-SM} - $\mathcal{S}'(\Rd)$-valued stable motion with finite dimensional distributions given by (\ref{def:lim_big}).\\
Let us notice first that the results for the case of finite and infinite variance are in a sense similar. The processes in Table 1 are Gaussian counterparts of stable processes in Table 2. Informally speaking, the finite variance branching law is a limit of laws given by (\ref{def:generating_beta}) hence one can observe similar phenomena in both cases. The case of intermediate dimensions is most interesting. The limits have similar spatial structure and complicated temporal one with long-range dependence property. The dependence on the starting distribution is intriguing since $Eq$ measure is the limit for a Poisson-starting system (\ref{def:Eq}) and by the time did not acquire any intuitive explanation. However both limits have the long-range dependence property  but for the equilibrium-starting system this dependence is stronger.\\
The remarkable feature of the limit process in the equilibrium case is that it can be decomposed into a sum of two independent $(\pb)$-stable processes. One of them being exactly the limit process in the Poisson case. This decomposition is an analogue of the one studied in \cite{Dzhaparidze:2004aa} for fractional Brownian motion. It should be noted that the process obtained in this case is a stable analogue of fractional Brownian motion. Namely, it is self-similar and has stationary increments. Processes with this properties were discussed in \cite[Chapter 7]{Samorodnitsky:1994aa} (see also Remark \ref{rem:Hssi}). This fact makes analogies between infinite-variance and finite-variance cases even stronger (recall that the limit process in the case of the finite variance branching law is fractional Brownian motion).
 The cases of critical and large dimensions are less complicated. The limits have independent increments and complicated spatial structure for large dimensions going beyond the space of measures. The qualitative change of the type of limits with dimension can be, partially, explained by recurrence and transient property of the underlying $\alpha$-stable L\'evy motion.
Further results on the fluctuations of the occupation time can be found in \cite{bojdecki-2006}, \cite{bojdecki-2007-12}, \cite{Bojdecki:2007ab} were high-density limits and system with inhomogeneous starting distributions are studied. One should also mention \cite{Birkner:aa} and \cite{Birkner:2005aa} where similar problems are considered in discrete setting (lattice $\mathbf{Z}^d$).\\
Although the proofs in this paper  rely mostly on the schema and methods used in \cite{Bojdecki:2007aa}, \cite{bojdecki-2005} and \cite{Milos:2007aa} we had to overcome some new technical difficulties which emerged during studies on additional terms arising in analysis of  equilibrium-starting system.
	
\section{Results} \label{sec:results} 
	By $\mathcal{T}$ we denote the semigroup of the $\alpha$-stable motion and by $p_t$ its transition density, i.e.,
	\begin{equation}
		\mathcal{T}_t f(x) = (p_t \ast f)(x). \label{def:T}
	\end{equation}
	Let $M$ be an independently scattered random $(\pb)$-stable measure on $\Rd \times \mathbb{R}_+$ with the Lebesgue control measure. More precisely, for a Borel set $A$, $M(A)$ is $(\pb)$-stable variable with characteristic function
	\[ \exp \left\lbrace \lambda(A) |z|^{\pb} \left( 1 - i \sgn(z) \tan \frac{\pi}{2}(\pb)\right)\right\rbrace, \]
	variables on disjoint set are independent and $M$ is $\sigma$-additive a.s.\\
	We define two stable processes.
	\begin{equation}
	  \eta^{1}_{t} = \int_{\mathbb{R}^{d+1}} \rbr{\mathbf{1}_{[0,t]}(r) \int_{r}^{t}p_{u-r}(x)du} M(dx, dr) \label{def:eta1}
	\end{equation}
	and
	\begin{equation}
	 \eta^{2}_{t} =\inti \intr \left( \intc{t} p_{s+l}(x) ds \right) M(dx,dl) \label{def:eta2}
	\end{equation}
	It can be checked that for intermediate dimensions both processes are well-defined in the sense given in \cite[Chapter 3]{Samorodnitsky:1994aa}.\\
	Assume now that $\eta^1$ and $\eta^2$ are independent, then we define 
	\begin{equation}
	 \eta := \eta^1 + \eta^2. \label{def:eta}
	\end{equation}
 This process plays fundamental role in this paper. Detailed presentation of its properties is postponed to Theorem \ref{thm:eta_hssi} and Proposition \ref{prep:eta_kappa}. Now we give series of three theorems which are the main results of the paper.\\
In these theorems  $X_{T}$ is the rescaled occupation time fluctuation process defined by (\ref{def:occupation_fluct}) for a system $N$ starting from equilibrium  distribution (\ref{eq: laplace_equilibrum}).
	\begin{thm} \label{thm:intermediate}
Assume $\frac{\alpha}{\beta} < d < \frac{\alpha (\pb)}{\beta}$ and  $F_{T} = T^{(2+\beta - \frac{d}{\alpha}\beta)/(\pb)}$. Then
		\[ X_{T} \Rightarrow_c K\eta\lambda , \]
		where
		\[
			K = \rbr{-\frac{V}{\pb}\cos \frac{\pi}{2}(\pb)}^\pb.
		\]  
	\end{thm}
\begin{rem}
As announced in the Introduction the process $\eta$ consists of two independent summands (\ref{def:eta}) where $\eta^1$ is the process that occurs in the limit for the Poisson system (see \cite{Bojdecki:2007aa}). An explanation of the reason of this structure of $\eta$ as well as the interpretation of $\eta^2$ require further studies. If $\beta=1$ and $\eta$ is a fractional Brownian motion on the the whole line then $\eta^1 = \frac{\eta_t + \eta_{-t}}{2}$ (sub-fractional Brownian motion) and $\eta^1 = \frac{\eta_t - \eta_{-t}}{2}$ (see \cite{Dzhaparidze:2004aa}).
\end{rem}

	\begin{thm}\label{thm:critical}
	 Assume $d = \frac{\alpha (\pb)}{\beta}$ and $F_{T} = (T \log T)^{\frac{1}{\pb}}$. Then
		\[ X_T \rightarrow_i K \lambda \xi \text{ and } X_T \rightarrow_f K \lambda \xi \text{ as } T \rightarrow +\infty \]
		where $\xi$ is $(\pb)$-stable process with stationary independent increments and characteristic function
		\begin{equation}
		 \ev{\exp(iz\xi_t)} = \exp\left\lbrace -t|z|^{\pb} \left( 1-i \sgn(z) \tan\frac{\pi}{2}(\pb)\right)\right\rbrace, z\in \mathbf{R}, t\geq 0 \label{def:lim_critical}
		\end{equation}
		and
		\[ K = \left(-V \cos\frac{\pi}{2}(\pb)\intr \left( \int_0^1 p_r(x) dr\right)^\beta p_1(x) dx \right)^{\frac{1}{\pb}} \]
	\end{thm}
Before presenting the last theorem we introduce the potential operator corresponding to the $\alpha$-stable motion.
	\[ \mathcal{G}f(x) = \int_0^{+\infty}\mathcal{T}_t f(x) dt. \]
	\begin{thm}\label{thm:big}
Assume $d > \frac{\alpha (\pb)}{\beta}$ and $F_{T} = T^{\frac{1}{\pb}}$. Then
		\[ X_T \rightarrow_i X \text{ and } X_T \rightarrow_f X \text{ as } T \rightarrow +\infty, \]
		where $X$ is an $\mathcal{S}'(\Rd)$-valued $(\pb)$-stable process with stationary independent increments and characteristic function
		\begin{eqnarray}
			\ev{\exp(i\ddp{X(t)}{\phi})} = \exp\left\lbrace -K^{\pb}t\intr |\mathcal{G}\phi(x)|^{\pb} \left(1 - i(\sgn \mathcal{G}\phi(x) \tan \frac{\pi}{2}(\pb)) \right) dx \right \rbrace, \nonumber \\
		\phi \in \mathcal{S}(\Rd),t\geq 0, \: \label{def:lim_big}
		\end{eqnarray}
		where
		\[ K = \left( - \frac{V}{\pb} cos\frac{\pi}{2}(\pb) \right)^\pb \]
	\end{thm}

	\begin{rem}
		The limits in the last two theorems have independent increments and are non-Gaussian hence by \cite[Theorem 13.4]{Kallenberg:2002aa} are not continuous. It is somehow unexpected since processes $X_T$ are clearly continuous. This is also the reason why we can not obtain functional convergence in those cases.
	\end{rem}
Propositions below summarize basic properties of $\eta$ defined by (\ref{def:eta}).
		\begin{thm} \label{thm:eta_hssi}
			$\eta$ is a continuous, $(\pb)$-stable process process. It is self-similar with exponent $H = (2+\beta - \frac{d \beta}{\alpha} )/(1+\beta)$ and has stationary increments. 
		\end{thm}

The self-similarity can be proved by a simple calculation using the characteristic function of the  finite-dimensional distributions of (\ref{def:eta}) obtained using \cite[(3.2.2)]{Samorodnitsky:1994aa}. Processes $X_T$ have stationary increments which comes straightforward from the fact that $N$ is stationary (since it is a Markov process starting from stationary distribution). From Theorem \ref{thm:intermediate} we know that the process $\eta$ is a limit of $X_T$ hence is also stationary.

\begin{rem} \label{rem:Hssi}
		Processes of this type were discussed in \cite[Chapter 7]{Samorodnitsky:1994aa}. In the notation used there $\eta$ is $H$-ssi stable process. Contrary to the Gaussian case, where there is a unique up to a constant $H$-ssi process for a given $H$ (fractional Brownian motion with Hurst parameter $H$), there are plenty of stable $H$-ssi's. It would be interesting to check if $\eta$ is one of already known processes (this could draw analogies to other problems) or is a new process. Unfortunately, we do not know the answer to this question.
\end{rem}
	We introduce now a general notation to investigate long-range dependence in the case of stable processes (see \cite{Bojdecki:2007aa}).
	\begin{definition} \label{def:dep-exp}
		Let $\eta$ be a real infinitely divisible process. For $0\leq u< v<s<t, T>0, z_1,z_2 \in \mathbf{R}$ define
		\begin{eqnarray}
			D_T(z_1,z_2;u,v,s,t)&=& |\log \ev{e^{iz_1(\eta_v - \eta_u) + iz_2(\eta_{T+t} - \eta_{T+s})}}  \nonumber \\ 
			 &&- \log\ev{e^{iz_1(\eta_v - \eta_u)}} - \ev{e^{iz_2(\eta_{T+t} - \eta_{T+s})}} |,
		\end{eqnarray}
		Dependence exponent  $\kappa$ is defined by
		\[ \kappa = \inf_{z_1, z_2 \in \mathbf{R}} \inf_{0\leq u<v<s<t} \sup\lbrace \gamma>0: D_T(z_1,z_2;u,v,s,t) = o(T^{-\gamma}) \text{ as } T\rightarrow +\infty \rbrace \]
	\end{definition}
	By \cite[Theorem 2.7]{Bojdecki:2007aa} exponent $\tilde{\kappa}$ of $\eta_1$ (denoted by $\xi$ therein) is 
	\[
		\tilde{\kappa} = \left \lbrace 
		\begin{array}{cc}
			 \frac{d}{\alpha} & \beta > \frac{d}{d+\alpha}\\
			 \frac{d}{\alpha}\left( 1+\beta - \frac{d}{\alpha+d}\right) & \beta \leq \frac{d}{d+\alpha}.
		\end{array} 
		\right . 
	\]
Conducting similar computations as in \cite[Proof of Theorem 2.7]{Bojdecki:2007aa} it can be checked that dependence exponent for $\eta_2$ is $\kappa = \frac{d}{\alpha}-1$. 
It is straightforward consequence of Definition \ref{def:dep-exp} that the dependence exponent of a sum of independent processes is minimum of the exponents of the summands hence we obtain
	\begin{prep} \label{prep:eta_kappa}
		Process $\eta$ has dependence exponent $\kappa = \frac{d}{\alpha}-1$. 
	\end{prep}
	\begin{rem}
		It is interesting to notice that addition of an independent term $\eta_2$ arising in the limit for the equilibrium-starting system (recall Theorem \ref{thm:intermediate}) increases long-range dependence and the dependence exponent does not depend on $\beta$ any more.
	\end{rem}
	\section{Proofs}
	For the sake of brevity proofs for Theorem \ref{thm:critical} and \ref{thm:big} are omitted. They are direct combination of the methods of \cite{bojdecki-2005} and the argument employed in the proof of Theorem \ref{thm:intermediate}. The scheme below is quite general and could be easily adapted for those proofs.
	\subsection{Scheme of the proof} \label{sec:proof_schema}
	To make the proof clearer we present a general scheme. Detailed calculation are deferred to a  separated section. We treat measure-valued processes as $\mathcal{S}'(\Rd)$-valued one. This enables usage of space-time method from \cite{Bojdecki:1986aa}. Let $\tilde{X}_T$ denote $\mathcal{S}'(\Rd)$-random variable defined by (\ref{def: space-time}) corresponding to the process $X_T$. In order to prove weak convergence in $\mathcal{C}([0,1], \mathcal{S}'(\Rd))$ to $X$ it suffices to prove weak convergence of  $\tilde{X}_T$ 
	\begin{equation}
		\tilde{X}_T \Rightarrow \tilde{X} \label{lim:main2}
	\end{equation}
	and tightness of $\lbrace X_T \rbrace_{T\geq 1}$. In order to obtain (\ref{lim:main2}) it suffices to verify that
	\begin{equation}
		\ev{e^{\ddp{\tilde{X}_T}{\Phi}}} \rightarrow \ev{e^{\ddp{\tilde{X}}{\Phi}}} \label{lim:laplace}
	\end{equation}
	for any non-negative $\Phi\in \mathcal{S}(\Rd)$. The tightness can be proven utilizing the Mitoma theorem \cite{Mitoma:1983aa}, which states that tightness of $\lbrace X_T\rbrace_{T \geq 1} $ in $\mathcal{C}([0,\tau], \mathcal{S}'(\Rd))$ is equivalent to tightness of $\ddp{X_T}{\phi}$ in $\mathcal{C}([0,\tau], \Rd)$ for any $\phi \in \mathcal{S}(\Rd)$. 
	\subsubsection{Space-time convergence} \label{sec:space-time}
	 The purpose of this subsection is a calculation of the Laplace transform and gathering facts used to show convergence (\ref{lim:laplace}). The schema described below generally follows the lines of a scheme presented in \cite{Bojdecki:2007aa}, \cite{bojdecki-2005} and \cite{Milos:2007ab} hence we omit some details.\\
	To make the proof shorter we will consider $\Phi$ of the special form:
	\[ \Phi(x,t) = \varphi(x) \psi(t)\:\:\varphi\in \mathcal{S}(\Rd), \psi\in \mathcal{S}(\mathbb{R}^+), \varphi\geq 0, \phi \geq 0. \]
	We also denote 
	\begin{equation}
		\varphi_T = \frac{1}{F_T}\varphi, \: \chi(t) = \int_t^1 \psi(s) ds, \: \chi_T(t) = \chi(\frac{t}{T}). \label{def: abbrev} 
	\end{equation}
	We write  
	\begin{equation}
		\Psi(x,t) = \varphi(x) \chi(t), \label{def: Psi_bez_T} 
	\end{equation}
	\begin{equation}
		\Psi_T(x,t) = \varphi_T(x) \chi_T(t), \label{def: Psi} 
	\end{equation}
	note that $\Psi$ and $\Psi_T$ are positive functions. For generating function $F$ we define $G(s) = F(1-s) - (1-s)$ so in our case
	\[ G(s) = \frac{s^{\pb}}{\pb} \]
	Behavior of the system starting off from a single particle at $x$ is described by the function 
	\begin{equation}
		v_{\Psi}\left(x,r,t\right)=1-\mathbb{E}\exp\left\{ -\int_{0}^{t}\left\langle N_{s}^{x},\Psi\left(\cdot,r+s\right)\right\rangle ds\right\} \label{def: v} 
	\end{equation}
	where $N_{s}^{x}$ denotes the empirical measure of the particle system with the initial condition $N_{0}^{x}=\delta_{x}$. $v_{\Psi}$ satisfies the equation
	\begin{equation}
		v_{\Psi}\left(x,r,t\right)=\int_{0}^{t}\mathcal{T}_{t-s}\left[\Psi\left(\cdot,r+t-s\right)\left(1-v_{\Psi}\left(\cdot,r+t-s,s\right)\right)-VG\left(v_{\Psi}\left(\cdot,r+t-s,s\right)\right)\right]\left(x\right)ds.\label{eq: v} 
	\end{equation}
	This equation can be derived using the Feynman-Kac formula in the same way as \cite[Lemma 3.4]{Milos:2007ab}.
	 We also define 
	\begin{equation}
		n_{\Psi}\left(x,r,t\right)=\int_{0}^{t}\mathcal{T}_{t-s}\Psi\left(\cdot,r+t-s\right)\left(x\right)ds.\label{def: n} 
	\end{equation}
	Since we consider only positive $\Psi$, hence (\ref{def: v}) and (\ref{eq: v}) yield 
	\begin{equation}
		0\leq v_T(x,r,t) \leq n_T(x,r,t) \label{ineq: n<v}, 
	\end{equation}
	\begin{equation}
		v_T(x,r,t) \leq 1 \label{ineq: v<1},
	\end{equation}
	In the sequel, for simplicity of notation, we write 
	\begin{equation}
		v_{T}\left(x,r,t\right)=v_{\Psi_{T}}\left(x,r,t\right),\label{eq:v_T notacja} 
	\end{equation}
	\begin{equation}
		n_{T}\left(x,r,t\right)=n_{\Psi_{T}}\left(x,r,t\right),\label{eq:n_T notacja} 
	\end{equation}
	\begin{equation}
		v_{T}\left(x\right)=v_{T}\left(x,0,T\right),\label{eq: v_T_x notacja} 
	\end{equation}
	\begin{equation}
		n_{T}\left(x\right)=n_{T}\left(x,0,T\right)\label{eq: n_T_x notacja} 
	\end{equation}
	when no confusion arises.\\
	\begin{fact}
		$n_{T}\left(x,T-s,s\right)\rightarrow0$ uniformly in $x\in\mathbb{R}^{d}$, $s\in\left[0,T\right]$ as $T\rightarrow+\infty$. \label{fact: uniformaly} 
	\end{fact}
	The fact was proved in \cite[Fact 3.7]{Milos:2007ab}. From the proof therein we obtain also the inequality
	\begin{equation}
		n_T \leq \frac{c}{F_T}. \label{ineq:n_T<1/F_T} 
	\end{equation}
	
	Following the lines of \cite[Section 3.2.2]{Milos:2007ab} we introduce function $V_T$
	\begin{equation}
		V_T(x,l) = 1 - \ev{\exp\rbr{\ddp{N_l^x}{\ln (1-v_T)}}} \label{def:V} 
	\end{equation}
	which satisfies the equation (see \cite[(3.20)]{Milos:2007ab})
	\begin{equation}
		V_{T}\left(x,l\right)=\mathcal{T}_{l}v_{T}\left(x\right)-V\int_{0}^{l}\mathcal{T}_{l-s}G\left(V_{T}\left(\cdot,s\right)\right)\left(x\right)ds.\label{main_equation} 
	\end{equation}
	
	A trivial verification using (\ref{def:V}), (\ref{ineq: v<1}) and (\ref{main_equation}) provides us with
	
	\begin{equation}
		0\leq V_{T}\left(x,l\right)\leq\mathcal{T}_{l}v_{T}\left(x\right),\,\forall_{x\in\mathbb{R}^{d},l\geq0}. \label{ineq: VT<TvT} 
	\end{equation}
	Next we write the Laplace transform of the occupation time fluctuation process (\ref{def:occupation_fluct}) for the system $N$ starting from equilibrium distribution
	\begin{equation}
		\ev{e^{-\ddp{\tilde{X}_T}{\Phi}}} = e^{A(T) + B(T)} \label{eq:laplace-eq} 
	\end{equation}
	where 
	\begin{equation}
		A\left(T\right)= \int_{\mathbb{R}^{d}}\int_{0}^{T}\Psi_{T}\left(x,T-s\right)v_{T}\left(x,T-s,s\right)+VG\left(v_{T}\left(x,T-s,s\right)\right)dsdx, \label{def: A} 
	\end{equation}
	\begin{equation}
		B\left(T\right)= V\int_{0}^{+\infty}\int_{\mathbb{R}^{d}}G\left(V_{T}\left(x,t\right)\right)dxdt. \label{def: B} 
	\end{equation}

The derivation of this formula can be found in \cite[Section 3.2.2]{Milos:2007ab}.

	To show (\ref{lim:laplace}) we need to calculate limits of $A(T)$ and $B(T)$.  Let us notice here that $\exp(A(T))$ is the same as the right-hand side of \cite[(3.7)]{Bojdecki:2007aa} so it is the Laplace transform of a Poisson-starting system
	\[
		A(T) \rightarrow \frac{V}{\pb} \intr \intc{1} \left[\intr \int_r^1 \varphi(y) \psi(s) \int_r^s p_{u-r}(x) du ds dy\right]^{\pb} dr dx.
	\]
In Section \ref{sec:calculations} we shall prove
\begin{equation}
	B(T) \rightarrow \frac{V}{\pb} \inti \intr \left( \intc{1} \intr \varphi(y) p_{l+s}\left(x\right)\chi(s) ds dy \right)^\pb dx dl \label{eq:goal}.
\end{equation}
\\Finally, we are in position to interpret the result. By properties of the Laplace transform, the limit of (\ref{eq:laplace-eq}) splits into two independent parts, corresponding to $A(T)$ and $B(T)$ respectively. The first was investigated in \cite{Bojdecki:2007aa} and corresponds to the process $\eta^1$ defined by (\ref{def:eta1}).\\
Let us now investigate the second one. Denote the limit of (\ref{eq:goal}) by $B$. It can be handled in the following way
\[ B = \frac{V}{\pb} \inti \intr \left( \intr \intc{1} p_{l+s}(x) \int_s^1 \Phi(y, u) du ds dy\right)^ \pb dx dl. \]
We write
\[ F(x, l) = \intr \intc{1} p_{l+s}(x) \int_s^1 \Phi(y, u) du ds dy, \]
then
\[ B = \frac{V}{\pb} \inti \intr F(x,l)^\pb dx dl. \]
An argument as in \cite[Corollary 3.5]{Bojdecki:2007aa} implies that the characteristic function of $\ddp{\tilde{X}}{\Phi}$ is
\begin{equation}
	\exp\left\{-K \inti \intr |F(x,l)|^\pb \left(1-i\:\sgn (F(x,l))\tan\frac{\pi}{2}(\pb) \right)dx dl\right\}, \label{eq:tmp:char}
\end{equation}

where $K = \frac{V}{\pb} \cos\frac{\pi}{2}(\pb)$.\\ 
Following the lines of reasoning in \cite[End of Section 3]{Bojdecki:2007aa}, we can obtain the characteristic function of the  finite dimensional distributions (by passing to the limit with appropriate sequence approximating $ \Phi(y, s) = \sum_j z_j \varphi_j \delta_{t_j}(s)$). It is of the form (\ref{eq:tmp:char}) with 
\[ F(x,l) = \sum_j z_j \ddp{\lambda}{\varphi_j} \intc{t_j} p_{l+s}(x) ds. \]
Using theorem \cite[Proposition 3.4.2]{Samorodnitsky:1994aa} one can infer easily that $\eta^2$ (recall (\ref{def:eta2})) has the same finite dimensional distributions.

\subsubsection{Tightness} \label{sec:tightness_schema}
It has been already mentioned that to prove tightness it suffices to show tightness of real-valued processes $\ddp{X_T}{\varphi}$ for any $\varphi \in \mathcal{S}(\Rd)$. We apply below a scheme presented in \cite{Bojdecki:2007aa}. By \cite[Theorem 12.3]{Billingsley:1968aa} it is enough to show that there exist constants $\chi >0 $ and $\nu\geq 0$ such that
\begin{equation}
	\mathbb{P}(|\ddp{X_T(t_2)}{\varphi} - \ddp{X_T(t_1)}{\varphi}| \geq \delta) \leq \frac{C(\varphi)}{\delta^\nu}(t_2 - t_1)^{1+\chi} \label{ineq:bill}
\end{equation}
holds for all $t_1,t_2\in[0,1]$, $t_1<t_2$, all $T\geq1$, and all $\delta>0$. A lemma in \cite[Section 3]{Bojdecki:2006ab} shows that each $\varphi\in \mathcal{S}(\Rd)$ can be decomposed $\varphi = \varphi_1 - \varphi_2$, $\varphi_1,\varphi_2 \in \mathcal{S}(\Rd)$, and $\varphi_1,\varphi_2 \geq 0$, hence from now on we will assume that $\varphi \geq 0$. 
Tail probability can be estimated using inequality \cite[(3.39)]{Bojdecki:2006ab}
\begin{equation}
	\mathbb{P}(|\ddp{\tilde{X}_T}{\varphi \otimes \phi}| \geq \delta) \leq C \delta \intc{1/\delta}(1-Re(\ev{\exp(-i \theta \ddp{\tilde{X}_T}{\varphi \otimes \phi} )})) d \theta \label{ineq:char}
\end{equation}
Now an analysis similar to that in \cite{Bojdecki:2006aa} shows us (\ref{ineq:bill}). Indeed, we approximate $\delta_{t_2} - \delta_{t_1}$ by $\phi \in \mathcal{S}(\Rd)$ such that
$\chi(t) = \int_t^1 \phi(s) ds$ fulfils
\[
	0\leq \chi \leq \mathbf{1}_{[t_1,t_2]}.
\]
Notice that $\ddp{\tilde{X}_T}{\varphi \otimes \phi}|$ approximates $|\ddp{X_T(t_2)}{\varphi} - \ddp{X_T(t_1)}{\varphi}|$.\\
Suppose that we know that 
\begin{equation}
	\delta \intc{1/\delta}(1-Re(\ev{\exp(-i \theta \ddp{\tilde{X}_T}{\varphi \otimes \phi} )})) d \theta \leq \frac{C(\varphi)}{\delta^\nu}(t_2 - t_1)^{1+\chi} \label{ineq:tightness_main}
\end{equation}
A passage to the limit using Fatou's lemma implies (\ref{ineq:bill}). The task of proving the last inequality is delegated to Section (\ref{sec:tightness_calculations}). The proof there will be conducted by estimating characteristic function appearing on the right-hand side of (\ref{ineq:char}). Analogously to (\ref{eq:laplace-eq}) we have
\begin{equation}
 	\ev{\exp\left(-i \ddp{\tilde{X}_T}{\varphi \otimes \phi} \right)} = \exp\left( A(T) + B(T)\right), \label{eq:decomp_complex}
\end{equation}
where $A(T)$ and $B(T)$ are given by (\ref{def: A}) and (\ref{def: B}) with $v_T$ and $V_T$ being the complex counterparts of functions $v$ and $V$ from Section \ref{sec:space-time}. It is easy to check that they fulfil equations (\ref{def:V})
and
\[
	v_T(x,t) = \intc{t} \mathcal{T}_{t-s}\left[i \varphi_T(\cdot)\chi_T(T-s)(1-v_T(\cdot,s)) - G(v_T(\cdot,s))
	 \right](x)ds
\]
(cf. (\ref{eq: v})).

\subsubsection{Auxiliary facts}
Before proceeding to calculations we gather a few additional facts.
	Let $p_{t}$ denote transition density of $\alpha$-stable motion. We have 
	\begin{equation}
		p_t(x) = t^{-\frac{d}{\alpha}}p_1(xt^{-\frac{1}{\alpha}}) \label{eq:p_t p_1} 
	\end{equation}
	And hence 
	\begin{equation}
		\Vert p_t \Vert _q ^q = t^{\frac{d}{\alpha}(1-q)} \Vert p_1 \Vert _q ^q,\: q> \frac{d}{d+\alpha}. \label{eq:norm_p} 
	\end{equation}
	We will need a few straightforward inequalities 
	\begin{equation}
		(a+b)^\pb \leq 2^\pb(a^\pb + b^\pb), \label{ineq:a+b pb<} 
	\end{equation}
	\begin{equation}
		(a+b)^\pb - a^\pb - b^\pb \geq \beta b^\beta a,\: b\geq a\geq 0, \label{ineq: a+b >} 
	\end{equation}
	\begin{equation}
		(a+b)^\pb - a^\pb - b^\pb \leq (\pb)a^\delta b^{1+\beta-\delta}, \beta\leq \delta \leq 1, a,b\geq 0, \label{ineq:a+b<}
	\end{equation}
	Moreover, we will use the generalized Minkowski inequality
\begin{equation}
	\norm{\int \! f}{p}{} \leq \int\! \norm{f}{p}{}, p\geq 1\label{ineq:Minkowski}
\end{equation}
	and Young's inequality
\begin{equation}
	\norm{f \ast g}{q}{} \leq \norm{f}{p_1}{}\norm{g}{p_2}{}, \: \frac{1}{q} = \frac{1}{p_1} + \frac{1}{p_2} -1. \label{ineq:Young}
\end{equation}

\section{Calculations for the proof of Theorem \ref{thm:intermediate}} \label{sec:calculations}
The general schema presented in Section \ref{sec:proof_schema} left the main technical difficulty of the proof, namely convergence of (\ref{eq:goal}) untouched. We decompose $B(T)$ in the following way
\[ B(T) = \frac{V}{\pb} \left(B_3(T) - B_2(T) - B_1(T) \right), \]
where
\[ B_1(T) = \inti \intr (\T{l}v_T(x))^{1+\beta} - (V_T(x,l))^{1+\beta} dx dl, \]
\[ B_2(T) = \inti \intr (\T{l}n_T(x))^{1+\beta} - (v_T(x,l))^{1+\beta} dx dl, \]
\[ B_3(T) = \inti \intr (\T{l}n_T(x))^{1+\beta} dx dl, \]
(\ref{eq:goal}) will be shown once we obtain
\begin{equation}
	 B_1(T) \rightarrow 0 \label{lim_B1}, 
\end{equation}
\begin{equation}
	 B_2(T) \rightarrow 0, \label{lim_B2}
\end{equation}
\begin{equation}
	 B_3(T) \rightarrow \norm{\varphi}{}{1+\beta}\inti \intr \left( \intc{1} p_{l+s}\left(x\right)\chi(s) ds\right)^\pb dx dl. \label{lim_B3}
\end{equation}

The proof of (\ref{lim_B2}) is omitted since it is similar to the case of (\ref{lim_B1}) but simpler.

\subsection{Convergence of $B_3$} 
By definition of $n_{T}$ (see (\ref{def: n})) we obtain
\[ B_3(T) = \inti \intr \left(\T{l} \intc{T} \T{s}\varphi_T(x)\chi_T(s) ds\right)^\pb dx dl. \]
Changing integration variable $s \rightarrow Ts$ and using definition of $\varphi_T$ (recall (\ref{def: abbrev}))
\[ B_3(T) = T^{\frac{d}{\alpha}\beta - 1}\inti \intr \left( \intc{1} \T{l+Ts}\varphi(x)\chi(s) ds \right)^\pb dx dl.\]
Using the definition of semigroup $\mathcal{T}$ (see (\ref{def:T})) yields
\[ B_3(T) = T^{\frac{d}{\alpha}\beta - 1}\inti \intr \left( \intc{1} \intr p_{l+Ts}(x-y) \varphi (y)\chi(s) dy ds\right)^\pb dx dl.\]
By (\ref{eq:p_t p_1}) we have
\[ B_3(T) = T^{\frac{d}{\alpha}\beta - 1}\inti \intr \left( \intc{1} \intr T^{-\frac{d}{\alpha}} p_{l/T + s}\left(T^{-\frac{1}{\alpha}}(x-y)\right) \varphi(y)\chi(s) dy ds\right)^ \pb dx dl ,\]
and, after obvious substitutions 
\[ B_3(T) = \inti \intr \left( \intc{1} \intr p_{l+s}\left(x-T^{-\frac{1}{\alpha}}y\right) \varphi (y)\chi(s) dy ds\right)^\pb dx dl.\]
This can be written as 
\[
	B_3(T) =  \inti \intr \rbr{f_l \ast g_T}^\pb dx dl 
\]
where
$f_l(x) = \int_0^1 p_{s+l}(x) \chi(s) ds$ and $g_T(x) = T^{\frac{d}{\alpha}}\varphi(xT^{\frac{1}{\alpha}})$. \\
Firstly, using the Jensen inequality we easily  check
\[
	\norm{f_l}{\pb}{\pb} = \norm{\intc{1} p_{l+s}(x)}{\pb}{\pb} \leq \intr \intc{1}(p_{l+s}(x))^\pb dx ds = \intc{1} \norm{p_{l+s}}{\pb}{\pb} ds = 
\]
(by (\ref{eq:norm_p}))
\[
	c \intc{1} (l+s)^{-\frac{d}{\alpha}\beta} ds \leq c_1 l^{-\frac{d}{\alpha}\beta}.
\]
Secondly, using (\ref{ineq:Minkowski}) and (\ref{eq:norm_p}) we get
\[
	\norm{f_l}{\pb}{\pb} = \norm{\intc{1}\! p_{l+s}(x)}{\pb}{\pb} \leq \rbr{\intc{1} \!\norm{p_{l+s}}{\pb}{} ds}^\pb \!\!\! =  \rbr{\intc{1} (l+s)^{-\frac{d}{\alpha} \frac{\beta}{\pb}} ds}^\pb \!\!\!\leq c
\]
Combining the last two estimates we get $\norm{f_l}{\pb}{\pb} \leq c (1\wedge l^{-\frac{d}{\alpha}\beta})$.
In this way we have proved that $f_l$ is $(\pb)$-integrable with respect to $x$ and $l$ since $\frac{d}{\alpha}\beta > 1 $.

Taking into account the form of $g_T$ (informally speaking $g_T$ converges to $\delta_x \cdot\norm{\varphi}{}{}$) we acquire the $\mathcal{L}^{\pb}$ convergence
\[
f_l\ast g_T \rightarrow f_l\cdot\norm{\varphi}{}{}
\]
which is exactly (\ref{lim_B3}).
\subsection{Convergence of $B_1$} 
We prove (\ref{lim_B1}). Applying (\ref{ineq:a+b<}) we write
\[ B_1(T) \leq B_{11}(T) + B_{12}(T), \]
where
\[ B_{11}(T) = \inti \intr (1+\beta) [ \T{l}v_T(x) - V_T(x,l)]^{\frac{1+\beta}{2}} [\T{l}v_T(x)]^{\frac{1+\beta}{2}} dx dl, \]
\[ B_{12}(T) = \inti \intr [\T{l}v_T(x) - V_T(x,l)]^\pb dx dl. \]
Using  (\ref{main_equation}), (\ref{ineq: VT<TvT}) and (\ref{ineq: n<v}) we get
\[ B_{12}(T) \leq c \inti \intr \rbr{\intc{l} \T{l-s} (\T{s} n_T(x))^\pb ds}^\pb dx dl \]
By definition of $n_T$ (see (\ref{def: n}))
\[ B_{12}(T) \leq c \inti \intr \rbr{\intc{l} \T{l-s} \rbr{\T{s} \intc{T} \T{u} \varphi_T(x) du}^\pb ds}^\pb dx dl \]
Combining with definition of $\mathcal{T}$ (\ref{def:T}) and subsituting $u \rightarrow Tu$ we can rewrite
\[ T^{(1+\beta)(\frac{d}{\alpha}\beta-1)} \inti \intr \rbr{\intc{l} \intr p_{l-s}(x-y) \rbr{\intc {1} \intr p_{Tu+s}(y-z) \varphi(z) dz du}^\pb dy ds}^\pb dx dl \]
Application of  (\ref{eq:p_t p_1}) and substitutions $y \rightarrow T^{-\frac{1}{\alpha}}y$ and  $x \rightarrow T^{\frac{1}{\alpha}}x$  yield
\[ T^{C_T} \inti \intr \rbr{\intc{l} \intr p_{l-s}\rbr{T^{\frac{1}{\alpha}}(x-y)} \rbr {\intc{1} \intr p_{u+s/T}\rbr{y-T^{-\frac{1}{\alpha}}z} \varphi(z) dz du}^\pb dy ds}^\pb dx dl \]
where $C_T=-(\pb)^2+\frac{d}{\alpha}(1+\beta^2 + \beta^3)$. Applying (\ref{eq:p_t p_1}) and substituting $s \rightarrow Ts$ and $l \rightarrow Tl$ we obtain
\[ T^{1-\frac{d}{\alpha}\beta} \inti \intr \rbr{\intc{l} \intr p_{l-s}\rbr{x-y} \rbr{\intc{1} \intr p_{u+s}\rbr {y-T^{-\frac{1}{\alpha}}z} \varphi(z) dz du}^\pb dy ds}^\pb dx dl \]
Let us denote
\[ f_T(s, y) = \rbr{\intc{1} \intr p_{u+s}\rbr{y-T^{-\frac{1}{\alpha}}z} \varphi(z) dz du}^\pb =  \]
\[  \rbr{\Phi_T \ast \intc{1} p_{u+s} du}^\pb (y), \]
where $\Phi_T(x) = T^{\frac{d}{\alpha}}\varphi(T^{\frac{1}{\alpha}}x)$. We obtain
\begin{equation}
B_{12}(T) \leq C_T \inti H(l) dl, \label{integral_to_be_proved}	
\end{equation}
where
\[ H(l) = \norm{\intc{l} p_{l-s} \ast f_T(s, \cdot) ds}{\pb}{\pb}. \]
Applying (\ref{ineq:Minkowski}) we get
\[ H(l) \leq \rbr{\intc{l} \norm{p_{l-s} \ast f_T(s, \cdot)}{\pb}{}ds}^\pb. \]
Utilizing Youngs's inequality (\ref{ineq:Young}) we write
\[ H(l) \leq \rbr{\intc{l} \norm{p_{l-s}}{\pb}{} \norm{f_T(s, \cdot)}{1}{} ds}^\pb. \]
By (\ref{eq:norm_p}) and the definition of $f_T$ we obtain
\[ H(l) \leq \rbr{\intc{l} (l-s)^{-\frac{d\beta}{\alpha (\pb)}} \norm{\Phi_T \ast \intc {1} p_{u+s} du}{\pb}{\pb} ds}^\pb. \]
Using (\ref{ineq:Young}) once again we have
\[ H(l) \leq \rbr{\intc{l} (l-s)^{-\frac{d\beta}{\alpha (\pb)}} \norm{\Phi_T}{1}{\pb} \norm{\intc{1} p_{u+s} du}{\pb}{\pb} ds}^\pb. \]
Hence by the (\ref{ineq:Minkowski}) and (\ref{eq:norm_p})
\[ H(l) \leq c \rbr{\intc{l} (l-s)^{-\frac{d\beta}{\alpha (\pb)}} \rbr{\intc{1} (u+s)^{- \frac{d\beta}{ \alpha (\pb)}} du}^\pb ds}^\pb. \]
For intermediate dimensions $\frac{d}{\alpha}\frac{\beta}{\pb}<1$, hence the inner integral can be estimated by a constant independent of $s$ so
\[ H(l) \leq c \rbr{\intc{l} (l-s)^{-\frac{d\beta}{\alpha (\pb)}} ds}^\pb. \]
Assume $l\leq 1$. Then
\begin{equation}
H(l)<c.	\label{H_small}
\end{equation}
Now we derive estimation that works for "large" $l$'s. By (\ref{ineq:Minkowski}) we have
\[ H(l) \leq \rbr{\intc{l} \norm{p_{l-s} \ast f_T(\cdot, s)}{\pb}{} ds}^\pb. \]
Young's inequality (\ref{ineq:Young}) yields
\[ H(l) \leq \rbr{\intc{l} \norm{p_{l-s}}{\pb}{} \norm{f_T(\cdot, s)}{1}{} ds}^\pb. \]
We estimate $\mathcal{L}^1$ norm of $f_t$ using (\ref{ineq:Young})
\[ \norm{f_T(\cdot, s)}{1}{} = \norm{\Phi_T \ast \intc{1} p_{u+s} du}{\pb}{\pb} \leq \norm{\Phi_T}{1}{\pb} \norm{\intc{1} p_{u+s} du}{\pb}{\pb}.\]
We use the trivial fact that $\norm{\Phi_T}{1}{}=c$ and (\ref{ineq:Young}) ($p,q\in [1,\pb]$ are yet to be specified)
\[ \norm{f_T(\cdot, s)}{1}{} \leq c \norm{p_s \ast \intc{1} p_{u} du}{\pb}{\pb} \leq \norm{p_s}{p}{\pb} \norm{\intc{1} p_{u} du}{q}{\pb}.\]
We have
\[ \norm{\intc{1} p_{u} du}{q}{\pb} \leq \rbr{\intc{1} \norm{p_{u}}{q}{} du}^ {\pb} \leq \rbr{\intc{1} u^{\frac{d}{\alpha}(\frac{1}{q}-1)} du}^{\pb} \leq c, \]
since it is obvious that for any $q\in [1,\pb]$ we have $\frac{d}{\alpha}(\frac{1}{q}-1)>-1$. And finally
\[  \norm{f_T(\cdot, s)}{1}{} \leq c s^{\frac{d}{\alpha}\rbr{\frac{1}{p} - 1}\rbr{\pb}}. \]
One can adjust $p$ to make exponent $A={\frac{d}{\alpha}\rbr{\frac{1}{p} - 1} \rbr{\pb}}$ arbitrary near $-1$ (because if $p=\pb$ then $A<-1$). Hence going back to estimation of $H(l)$
\[ H(l) \leq c\rbr{\intc{l} (l-s)^B s^A ds}^\pb, \]
where $B=\frac{d}{\alpha}\rbr{\frac{1}{\pb}-1}$. Substituting $s\rightarrow ls$ we obtain
\[ H(l) \leq c \rbr{l^{A+B+1}\intc{1} (1-s)^B s^A ds}^\pb = c l^{(A+B+1)(\pb)}\rbr{\intc{1} (1-s)^B s^A ds}^\pb \]
Notice that $B>-1$ and
\[(A+B+1)(\pb) = (A+1)(\pb) + \rbr{-\frac{d}{\alpha}\beta} \]
$-\frac{d}{\alpha}\beta<-1$ and $A+1$ can be made arbitrarily near $0$, hence
\begin{equation}
H(l) \leq l^W	\label{H_big}
\end{equation}
where $W<-1$. Combining estimates (\ref{H_small}) and (\ref{H_big}) for $H$ we can conclude that the integral in (\ref{integral_to_be_proved}) is finite and $B_{12}\rightarrow 0$.

\subsection{Tightness calculations} \label{sec:tightness_calculations}
Following the scheme in Section \ref{sec:tightness_schema} it remains to prove (\ref{ineq:tightness_main}).
Slightly abusing notation we will use an additional argument $\theta$ to indicate that a function is computed for $\theta \Phi$ instead of $\Phi$ (eg. $A(T,\theta)$). In this section we deal with complex functions defined in Section \ref{sec:tightness_schema}, which are not to be confused with functions in sections devoted to space-time convergence. \\
From equation (\ref{ineq:tightness_main}) we have to estimate
\begin{equation*}
L = 1-Re\rbr{\ev{\exp\rbr{-i \theta \ddp{\tilde{X}_T}{\varphi \otimes \phi} }}} 
\end{equation*}
We use (\ref{eq:decomp_complex}) and then estimate $A(T,\theta)$ in the same way as in the proof of tightness in \cite{Bojdecki:2007aa}, thus
\[L \leq \ev{\left| 1-\exp\rbr{-i \theta \ddp{\tilde{X}_T}{\varphi \otimes \phi} } \right|} \leq |I + II + III| \]
where
\[ I = i \theta \intr \intc{T} \varphi_T(x) \chi_T(T-s) v_{T,\theta}(x,s) dx ds \]
\[ II = \frac{V}{\pb} \intr \intc{T} v_{T,\theta}^\pb(x,s) dx ds \]
\[ III = V \inti \intr \rbr{V_{T,\theta}(x,t)}^\pb dx dt \]
The terms  $I$ and $II$ are the same as in  \cite{Bojdecki:2007aa}. Hence we have only to deal with $|III|$. Before that we  show an estimation (which holds for $T$ large enough). Firstly, recall the definition of $V_T$ (\ref{def:V})
\[ |V_{T,\theta}(x, t)| = |1 - \ev{e^{\ddp{N^x_t}{ \ln w_{T,\theta} }}}| \leq \ev{|1 - e^{\ddp{N^x_t} { \ln w_{T,\theta} }}|}  \]
where
\begin{equation*}
 w_{T,\theta}(x,r,t) = \mathbb{E}\exp\left\{ -i\theta\int_{0}^{t}\left\langle N_{s}^{x},\Psi\left(\cdot,r+s\right)\right\rangle ds\right\}.
\end{equation*}
We know that $|w_{T,\theta}| \leq 1$ which implies $|\ln w_{T,\theta}| \leq 0$ and consequently $e^{\ddp{N^x_t}{ \ln w_{T,\theta} }} \leq 1$. Finally, if $|z|<1$ we can use inequality $|1-e^z|\leq 2 |z|$. Hence
\[ |V_{T,\theta}(x, t)| \leq 2\ev{|\ddp{N^x_t}{ \ln w_{T,\theta} }|} \leq 2\ev{\ddp{N^x_t}{ |\ln w_{T,\theta}| }} = 2\T{t} |\ln w_T| \leq 2\T{t} n_{T,\theta}. \]
(note that $n_{T, \theta}$ is a real function).\\
Therefore we have to estimate
\[ |B_{T,\theta}| \leq \inti \intr (\T{l} n_{T,\theta}(x))^\pb dx dl \]
Let us notice that the integral is the same as $B_3$ from in Section \ref{sec:calculations}. For $T $ large enough we have
\[ |B_{T,\theta}| \leq C \norm{\theta \varphi}{}{1+\beta}\inti \intr \left( \intc{1} p_ {l+s}\left(x\right)\chi(s) ds\right)^\pb dx dl \]
According to the argument in Section \ref{sec:tightness_schema} we choose $\chi \simeq 1_{[t_1, t_2]}$
\[ |B_{T,\theta}| \leq\inti \intr \left( \int_{t_1}^{t_2} p_{l+s}\left(x\right) ds\right)^\pb dx dl \]
Denote $\Delta := t_2 - t_1$
\[ |B_{T,\theta}| \leq \inti \intr \left( \intc{\Delta} p_{l+s}\left(x\right) ds\right)^\pb dx dl \]
After obvious substitutions and using (\ref{eq:p_t p_1}) we obtain
\[ |B_{T,\theta}| \leq \Delta^{(2+\beta)} \Delta^{-\frac{d}{\alpha}\beta} \inti \intr \left( \intc{1} p_{l + s}\left( x\right) ds\right)^\pb dx dl. \]
It is easy to check that $2+\beta-\frac{d}{\alpha}\beta>1$ reasoning along the lines of the proof in \cite{Bojdecki:2007aa} completes the proof.\\

\bibliographystyle{amsplain}
\bibliography{branching}

\end{document}